# Differential algebras of Legendrian links[1]

Yu. V. Chekanov

June 1997

## 1 Introduction

Consider the space $\mathbf{R}^3 = \{(p,q,u)\}$ equipped with the standard contact form $\alpha = du - pdq$. A Legendrian link in $\mathbf{R}^3$ is a link $L$ such that the restriction of $\alpha$ to $L$ vanishes.

Consider the projection $\pi \colon \mathbf{R}^3 \to \mathbf{R}^2$, $(p,q,u) \mapsto (p,q)$. The projection $\pi(L)$ of a Legendrian link $L$ is an immersedcurve since the kernel of $dp$ is everywhere transverse to the kernel of $\alpha$. We say that $L$ is generic if all self-intersections of $\pi(L)$ are transverse double points. The diagram of a Legendrian link $L$ is its projection $\pi(L)$ pictured in such a way that at every double point the branch with the larger value of $u$ is shown as the upper one.

An oriented Legendrian link $L$ has three classical invariants: its isotopy class in the space of smooth embeddings, and Maslov and Bennequin numbers. The Maslov number $m(L_i)$ of the component $L_i$ of $L$ is twice the rotation number of $\pi(L_i)$. The Bennequin number $\beta(L_i)$ of $L_i$ is the linking number between $L_i$ and $s(L_i)$, where $s$ is a small shift along the $u$-axis, and $\mathbf{R}^3$ is oriented with the form $\alpha \wedge d\alpha = -dp \wedge dq \wedge du$. In terms of the diagram, $\beta(L_i)$ is the number of double points of $\pi(L_i)$ counted with signs shown in Fig. 1. The Bennequin number does not depend on the orientation of $L_i$.

---

[1] preliminary version

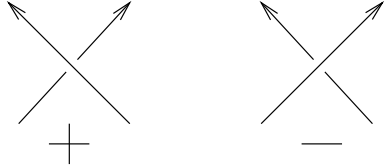

Figure 1:

The problem of classification of Legendrian knots (links) up to isotopy in the class of Legendrian embeddings (Legendrian isotopy) naturally leads to the following two subproblems. The first of them is: which combinations of the three classical invariants can be realized by a Legendrian knot? (It is well-known that each invariant by itself can be realized by a Legendrian knot). The first step in this direction was done by Bennequin who proved that the Bennequin number of a Legendrian knot is less than twice its genus [1] (see [2, 3, 4] for further results).

The second subproblem is whether there exist a pair of Legendrian knots which have the same classical invariants but are not Legendrian isotopic. Eliashberg and Fraser showed that this is impossible when knots are trivial as smooth knots [5, 6]. In the present paper, we develop a theory which allows, in particular, to prove the following assertion.

**Theorem 1.1** *Let $L_1, L_2 \subset \mathbf{R}^3$ be the Legendrian knots whose diagrams are represented in Fig. 6 on page 8. Then $m(L_1) = m(L_2)$, $\beta(L_1) = \beta(L_2)$ and $L_1$ is isotopic to $L_2$ in the class of smooth embeddings but not in the class of Legendrian embeddings.*

Our new invariants which tell those knots from each other appear in the following way. To every generic Legendrian link $L$ we associate a free differential algebra over $\mathbf{Z}_2$, whose generators are the double points of the diagram of $L$, and the differential is defined combinatorially, via the diagram. Legendrian isotopic links lead to in a sense equivalent differential algebras.

The differential algebra of a Legendrian link in $\mathbf{R}^3$ is a particular case of a more general Morse-theoretical construction which associates to every generic Legendrian submanifold $L$ in a contact manifold $(M, \alpha)$ (at least when the Reeb flow defined by $\alpha$ has no closed trajectories) a differential algebra whose generators are trajectories of the Reeb flow starting and ending on $L$. The differential is defined in terms of $J$-holomorphic curves in the symplectization of $(M, \alpha)$. It is known at least since Hofer's paper [7] on $J$-holomorphic curves in symplectizations that such a theory exists. There are no rigorous proofs at present, but recently Eliashberg and Hofer announced that they have obtained ones. On the combinatorial level, their theory is rather similar to ours. It also allows them to distinguish two knots (different from those of Theorem 1.1) whose classical invariants coincide.

The paper is organized as follows. In Section 1, we give the necessary algebraic definitions. In Section 2, we introduce the differential algebra of (a diagram of) a Legendrian link and formulate main theorems concerning it. In Section 3, a simple invariant of so-called stable tame isomorphism of free differential algebras is defined. Then it is applied to prove Theorem 1.1. Sections 4 and 5 contain the proofs of the theorems stated in Section 2. In Section 6, we discuss certain properties of our differential algebra.

**Acknowledgement** I am grateful to V. Goryunov and S. Tabachnikov for stimulating discussions and helpful remarks on the first draft of the paper. I also thank Ya. Eliashberg for stimulating discussions. It is a pleasure to acknowledge the hospitality of the Fields Institute for Research in Mathematical Sciences.

## 2 Free differential algebras and stable tame isomorphisms

Let $A$ be an associative $\mathbf{Z}_2$-algebra A pair $(A, d)$, where $d: A \to A$, is called a differential algebra if $d(ab) = d(a)b + a\,d(b)$ for every $a, b \in A$, and $d^2 = 0$.

Denote by $T(a_1, \ldots a_n)$ the free associative unitary algebra over $\mathbf{Z}_2$ with generators $a_1, \ldots a_n$, $n \geq 0$. There exists a vector space decomposition $T(a_1, \ldots a_n) = \oplus_{l=0}^{\infty} A_l$, where $A_l$ is spanned by monomials of the form $a_{i_1} \ldots a_{i_l}$ and $A_0$ is just $\mathbf{Z}_2$. If $(T(a_1, \ldots a_n), d))$ is a differential algebra, we call it a free differential algebra (the explicitly given generators $a_1, \ldots a_n$ here constitute an important part of the definition). Note that $d(1) = 0$ (and hence $d$ is linear), and $d$ is determined by its values on the generators $a_i$.

To every diagram $Y$ of a Legendrian link we associate a free differential algebra $(A_Y, d_Y)$ whose generators are the double points of $Y$. Since different diagrams can represent isotopic Legendrian links, we need to introduce an appropriate equivalence relation on the class of free differential algebras in order to construct an invariant of Legendrian isotopy.

First we define the stabilization operation. Given two differential algebras $(A, d_A) = (T(a_1, \ldots a_n), d_A)$, $(B, d_B) = (T(b_1, \ldots b_m), d_B)$, their coproduct is the differential algebra $(A, d_A) \amalg (B, d_B) = (T(a_1, \ldots a_n, b_1, \ldots b_m), d)$, where $d(a_i) = d_A(a_i)$, $d(b_j) = d_B(b_j)$. Consider the free differential algebra $(E, d_E)$, where $E = T(e_1, e_2)$, $d_E(e_1) = e_2$, $d_E(e_2) = 0$. We define the stabilization operation as follows: $S(A, d) = (A, d) \amalg (E, d_E)$. Denote by $S^N(A, d)$ the result of applying this operation $N$ times.

We call an automorphism $s$ of a free algebra $T(a_1, \ldots a_n)$ elementary if, for some $j \in \{1, \ldots n\}$, $s(a_i) = a_i$ when $i \neq j$, and $s(a_j) = a_j + u$, where $u \in T(a_1, \ldots a_{i-1}, a_{i+1}, \ldots a_n)$. The group of tame automorphisms of $T(a_1, \ldots a_n)$ is that generated by elementary automorphisms. (Actually, it is not known to the author whether non-tame automorphisms exist.)

We call an isomorphism $T(a_1, \ldots a_n) \to T(b_1, \ldots b_n)$ tame if it is a composition of a tame automorphism and the isomorphism sending $a_i$ to $b_i$ for every $i \in \{1, \ldots, n\}$. One easily checks that different orderings of the set of generators lead to tame isomorphic algebras.

A homomorphism (resp. isomorphism) $g: A \to B$ is called a homomorphism (resp. isomorphism) between the differential algebras $(A, d_A)$ and $(B, d_B)$ if $g \circ d_A = d_B \circ g$.

We call two differential algebras $(A, d) = (T(a_1, \ldots a_n), d)$, $(A', d') = (T(a'_1, \ldots a'_n), d')$ tame isomorphic if there exists a tame isomorphism between $T(a_1, \ldots a_n)$ and $T(a'_1, \ldots a'_n)$ which is an isomorphism of differential algebras. We call differential algebras $(T(a_1, \ldots a_n), d)$, $(T(a'_1, \ldots a'_{n'}), d')$ stable tame isomorphic if for some natural numbers $N, N'$ the stabilizations $S^N(T(a_1, \ldots a_n), d)$ and $S^{N'}(T(a'_1, \ldots a'_{n'}), d')$ are tame isomorphic.

## 3 Differential algebra of a diagram of Legendrian link

For every natural $k$, we fix a (curved) convex $k$-gon $\Pi_k \subset \mathbf{R}^2$ whose vertices $x_0^k, \ldots, x_{k-1}^k$ are numbered counterclockwise. Let $Y \subset \mathbf{R}^2$ be the diagram of a Legendrian link $L$; the symplectic form $dp \wedge dq$ equips $\mathbf{R}^2$ with an orientation. Denote by $W_k(Y)$ the collection of smooth orientation-preserving immersions $f\colon \Pi_k \to \mathbf{R}^2$ such that $f(\partial \Pi_k) \subset Y$ (and hence the images of the vertices $x_i^k$ are double points of $Y$). Let $G_k$ denote the subgroup of the group $\mathrm{Diff}_+(\Pi_k)$ of orientation-preserving diffeomorphisms of $\Pi_k$ consisting of the maps which fix the vertices $x_i^k$. Instead of the space $W_k(Y)$ itself, we shall use the set $\widetilde{W}_k(Y)$ of non-parametrized immersions $\widetilde{W}_k(Y) = W_k(Y)/G_k$.

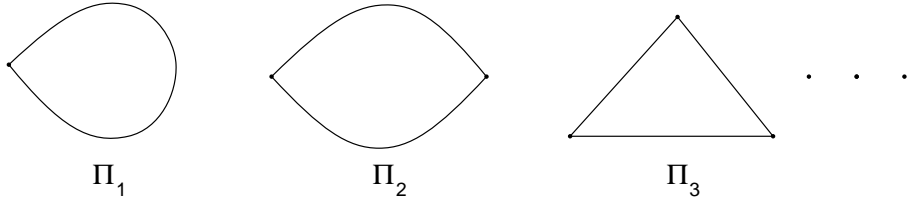

Figure 2:

Let $a_1, \ldots, a_n$ be the double points of $Y$. Consider the free algebra $A_Y = T(a_1, \ldots, a_n)$. We are going to equip the algebra $A_Y$ with a differential. For each of the double points $a_j$, the diagram $Y$ divides a neighbourhood of $a_j$ into four sectors. We mark two of them positive and the other two negative as shown below:

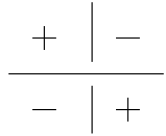

Figure 3:

The set $\widetilde{W}_k(Y)$ is a disjoint union of the sets $\widetilde{W}_k(Y, a_{j_0}, \ldots, a_{j_{k-1}})$ consisting of such $f$ that $f(x_i^k) = a_{j_i}$, $0 \leq i \leq k-1$. For each of the vertices $x_i^k$ of the polygon $\Pi_k$, the smooth immersion $f \in \widetilde{W}_k(Y)$ maps its neighbourhood to either positive or negative sector; we shall say that $x_i^k$ is correspondingly a positive or negative vertex for $f$.

We call an immersion $f \in \widetilde{W}_k(Y)$ admissible if the vertex $x_0^k$ is positive for $f$ and all other vertices are negative. Denote the set of admissible immersions by $W_k^+(Y)$; denote also $W_k^+(Y, a_{j_0}, \ldots, a_{j_{k-1}}) = W_k^+(Y) \cap$

$\widetilde{W}_k(Y, a_{j_0}, \ldots, a_{j_{k-1}})$. Define

$$d(a_j) = \sum_{k \geq 1} \sum_{f \in W_k^+(Y, a_{j_0}, \ldots, a_{j_{k-1}})} a_{j_1} \cdots a_{j_{k-1}} \qquad (1)$$

(the summands corresponding to $f \in W_1^+(Y, a_{j_0})$ are 1s) and extend $d$ to a linear map $A_Y \to A_Y$ satisfying the Leibnitz identity $d(ab) = d(a)b + a\,d(b)$. Note that the orientation of $L$ is not involved in this definition.

**Theorem 3.1** *The left hand side of (1) contains finitely many terms and hence the definition of $d$ is sound.*

**Theorem 3.2** $d^2 = 0$.

**Theorem 3.3** *Let $L_0, L_1$ be generic Legendrian links and let $(A_{Y_0}, d_{Y_0})$, $(A_{Y_1}, d_{Y_1})$ be the free differential algebras associated with their diagrams. If $L_0$ is Legendrian isotopic to $L_1$ then $(A_{Y_0}, d_{Y_0})$ is stable tame isomorphic to $(A_{Y_1}, d_{Y_1})$.*

The proofs of these assertions are given in Sections 5 and 6.

*Examples.* For several diagrams of Legendrian knots we compute their classical invariants and the associated differential algebras. Two of them represent the pair of knots of Theorem 1.1.

One can check that if one diagram of a Legendrian link can be mapped to another by an orientation-preserving diffeomorphism of $\mathbf{R}^2$ then the links are Legendrian isotopic. Thus we can unambiguously represent Legendrian links not only by their diagrams, but also by the images of the diagrams under orientation-preserving diffeomorphisms, without caring to reproduce the areas realistically.

1. The diagram $Y$ of the simplest Legendrian knot is represented in Fig. 4 The classical invariants of the corresponding Legendrian knot $L$ are as follows: $m(L) = 0$, $\beta(L) = -1$, the knot is the trivial knot in the smooth category. We have $A_Y = T(a_1)$. The set of immersions $\widetilde{W}_k(Y)$ is empty for $k > 1$ and consists of two elements $f_1, f_2 \in W_1^+(Y, a_1)$, whose images are the closures of the two bounded components of $\mathbf{R}^2 \backslash Y$. Hence $d_Y(a_1) = 1 + 1 = 0$.

2. Consider the diagram $Y_l$ given in Fig. 5a. As a smooth knot, the corresponding Legendrian knot $L_l$ is isotopic to the right-handed (with respect to the orientation given by $\alpha \wedge d\alpha = -dp \wedge dq \wedge du$) trefoil. We have $m(L_l) = 0$, $\beta(L_l) = 1$, $(A_{Y_l}, d_{Y_l}) = (T(a_1, \ldots, a_5), d)$, where

$$d(a_1) = a_4 + a_3 a_4 + a_4 a_5 + a_4 a_5 a_3 a_4$$
$$d(a_2) = 1 + a_3 + a_5 + a_3 a_4 a_5$$
$$d(a_3) = d(a_4) = d(a_5) = 0.$$

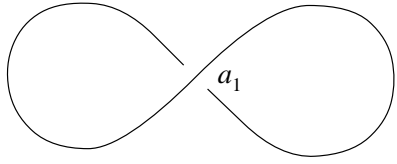

Figure 4:

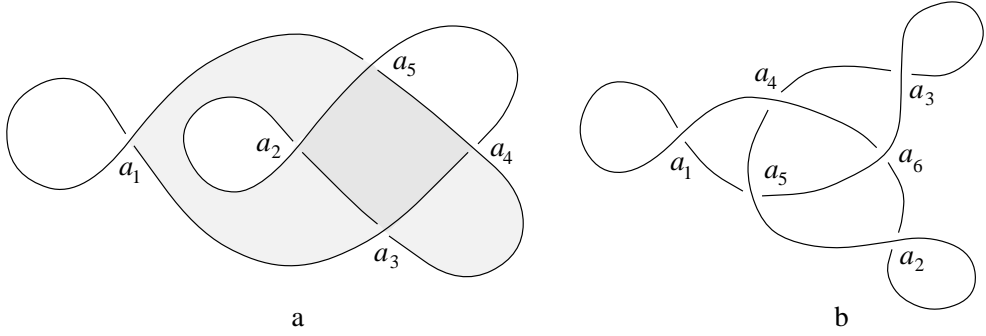

Figure 5:

For example, the shaded area is the image of the immersion corresponding to the term $a_4 a_5$ in $d_{Y_l}(a_1)$.

3. Consider the diagram $Y_r$ represented in Fig. 5b. The corresponding Legendrian knot $L_r$ is isotopic, in the class of smooth embeddings, to the left-handed trefoil. We have $m(L_r) = \pm 2$, $\beta(L) = -6$. $(A_{Y_r}, d_{Y_r}) = (T(a_1, \ldots, a_6), d)$, where

$$d(a_1) = 1 + a_4 a_6,$$
$$d(a_2) = 1 + a_5 a_4,$$

$$d(a_3) = 1 + a_6 a_5,$$
$$d(a_4) = d(a_5) = d(a_6) = 0.$$

4. Let $L_1$ be the Legendrian knot whose diagram $Y_1$ is given in Fig. 6a and $L_2$ the Legendrian knot whose diagram $Y_2$ is given in Fig. 6b. We have $m(L_1) = m(L_2) = 0$, $\beta(L_1) = \beta(L_1) = 1$. Clearly, the two knots are isotopic as smooth embeddings (they represent the class $5_2$ in the classification of smooth knots). Consider the differential algebra $(T(a_1, \ldots, a_9), d) = (A_{Y_2}, d_{Y_2})$, then

$$d(a_1) = 1 + a_7 + a_7 a_6 a_5,$$
$$d(a_2) = 1 + a_9 + a_5 a_6 a_9,$$
$$d(a_3) = 1 + a_8 a_7,$$
$$d(a_4) = 1 + a_8 a_9,$$
$$d(a_5) = d(a_6) = d(a_7) = d(a_8) = d(a_9) = 0.$$

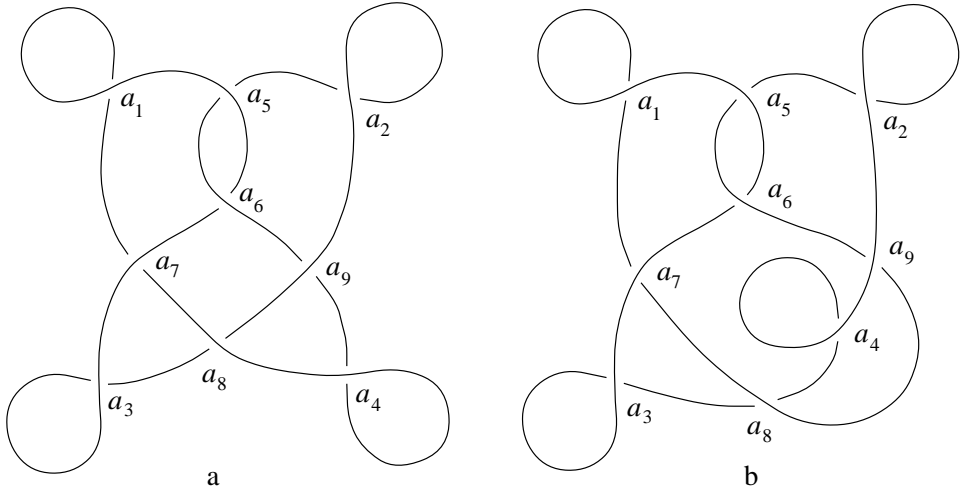

Figure 6:

For $(T(a_1, \ldots, a_9), d) = (A_{Y_2}, d_{Y_2})$, we have

$$d(a_1) = 1 + a_7 + a_7 a_6 a_5 + a_5,$$

$$d(a_2) = 1 + a_9 + a_5a_6a_9,$$
$$d(a_3) = 1 + a_8a_7,$$
$$d(a_4) = 1 + a_8a_9,$$
$$d(a_5) = d(a_6) = d(a_7) = d(a_8) = d(a_9) = 0.$$

We shall show in the next section that the differential algebras $(A_{Y_1}, d_{Y_1})$, $(A_{Y_2}, d_{Y_2})$ are not stable tame isomorphic and hence the Legendrian knots $L_1, L_2$ are not isotopic as Legendrian embeddings.

## 4  The invariant $I$

It is, in general, hard to find out whether or not two given free differential algebras are (tame) isomorphic. Hence there is a need for invariants allowing to distinguish non-isomorphic ones. The most obvious, and the most common, of such invariants is the cohomology $\ker d / \operatorname{im} d$ of a differential algebra $(A, d)$. However, in many cases the cohomology lacks a simple description. In the present section, we introduce an easy-to-compute invariant of stable tame isomorphism of free differential algebras.

Let $A = T(a_1, \ldots a_n)$, $\bar{A} = \oplus_{l=1}^{\infty} A_l$. The differential algebra $(A, d)$ is called augmented if $d(\bar{A}) \subset \bar{A}$. Let $d = \sum_{l=0}^{\infty} d_l$, where $d_l(a_i) \subset A_l$ for every $i \in \{1, \ldots, n\}$. Then $(A, d)$ is augmented if and only if $d_0 = 0$. We call a free differential algebra good if it is tame isomorphic to an augmented one, and bad otherwise.

*Remark.* Note that the cohomology of a good differential algebra does not vanish since 1 does not belong to the image of the differential.

Now we are ready to define, for a free differential algebra, the invariant $I$, which is a set of non-negative integers. Assume that the the differential algebra $(T(a_1, \ldots a_n), d)$ is augmented. Then $d(\bar{A}^m) \subset \bar{A}^m = \oplus_{i=m}^{\infty} A_i$ for every natural $m$. Hence $d$ gives rise to a linear operator $d_{(1)}$ on the quotient vector space $\bar{A}/\bar{A}^2$. Clearly, $d_{(1)}^2 = 0$. As long as the generators $a_1, \ldots, a_n$ are fixed, we can identify $\bar{A}/\bar{A}^2$ with the linear space $A_1$ spanned by $a_1, \ldots, a_n$. Then $d_{(1)}$ will coincide with the restriction of $d_1$ to $A_1$.

Let
$$i(T(a_1, \ldots a_n), d) = \dim(\ker d_{(1)}) - \dim(\operatorname{im} d_{(1)}) = n - 2\dim(\operatorname{im} d_{(1)}).$$

In other words, $i$ is the dimension of the cohomology of the operator $d_{(1)}$.

Denote by $\mathrm{Aut}_t(A)$ the group of tame automorphisms of $A = T(a_1, \ldots a_n)$. For an arbitrary free differential algebra $(T(a_1, \ldots a_n), d)$, we define $I(T(a_1, \ldots a_n), d)$ as a collection of the numbers $i(T(a_1, \ldots a_n), d^g)$, where $d^g = gdg^{-1}$, over all $g \in \mathrm{Aut}_t(A)$ such that $(T(a_1, \ldots a_n), d^g)$ is augmented. In particular, $I = \emptyset$ if and only if the differential algebra is bad.

**Lemma 4.1** *If free differential algebras $(A, d_A)$, $(B, d_B)$ are stable tame isomorphic then $I(A, d_A) = I(B, d_B)$.*

We can associate to every generic Legendrian link $L$ the set $I(L) = I(A_{Y(L)}, d_{Y(L)})$ of non-negative integers. Combined with Theorem 3.3, Lemma 4.1 claims that $I(L) = I(L')$ when $L$ is isotopic to $L'$ in the class of Legendrian embeddings. We call a Legendrian link $L$ good if $I(L) \neq \emptyset$.

Before proving Lemma 4.1, we are going to clarify the definition of $I$. Denote by $\mathrm{Aut}_t(A)$ the group of tame automorphisms of $A = T(a_1, \ldots a_n)$. Let the subgroup $\mathrm{Aut}_r(A) \subset \mathrm{Aut}_t(A)$ consist of such $g$ that $g(\bar{A}) \subset \bar{A}$. The subgroup $\mathrm{Aut}_0(A) \subset \mathrm{Aut}_t(A)$ consists of such $g$ that for every $i \in \{1, \ldots, n\}$ $g(a_i) = a_i + c_i$, where $c_i \in A_0$. As $\mathrm{Aut}_0(A)$ is isomorphic to $\mathbf{Z}_2^n$, the following assertion greatly simplifies the computation of $I$.

**Lemma 4.2** *The invariant $I(T(a_1, \ldots a_n), d)$ is a collection of the numbers $i(T(a_1, \ldots a_n), d^g)$ over all $g \in \mathrm{Aut}_0(A)$ such that $(T(a_1, \ldots a_n), d^g)$ is augmented.*

The key ingredient in the proof is the following

**Lemma 4.3** $\mathrm{Aut}_t(A) = \mathrm{Aut}_0(A) \mathrm{Aut}_r(A) = \mathrm{Aut}_r(A) \mathrm{Aut}_0(A)$.

*Proof.* Let $s \in \mathrm{Aut}_t(A)$ be an elementary automorphism which sends $a_i$ to $c_i + v_i$, where $c_i \in A_0$, and $v_i \in \bar{A}$ and fixes other generators Consider the elementary automorphisms $s' \in \mathrm{Aut}_0(A)$, $\bar{s} \in \mathrm{Aut}_r(A)$ which map $a_i$ to, respectively, $a_i + c_i$ and $v_i$ and fix other generators. Then $s = s' \circ \bar{s} = \bar{s} \circ s'$.

Let $g \in \mathrm{Aut}_t(A)$. We have

$$g = s_1 \circ \cdots \circ s_m = s'_1 \circ \bar{s}_1 \circ \cdots \circ s'_m \circ \bar{s}_m,$$

where $s_i \in \mathrm{Aut}_t(A)$, $s'_i \in \mathrm{Aut}_0(A)$, $\bar{s}_i \in \mathrm{Aut}_r(A)$, are elementary.

For any $g_0 \in \mathrm{Aut}_0(A)$ and any elementary $\bar{s} \in \mathrm{Aut}_r(A)$ we have $g_0 \circ \bar{s} \circ g_0^{-1} \in \mathrm{Aut}_0(A)$. Hence there exist $g_-, g_+ \in \mathrm{Aut}_0(A)$ such that

$$g = g_- \circ \bar{s}_1 \circ \cdots \circ \bar{s}_m = \bar{s}_1 \circ \cdots \circ \bar{s}_m \circ g_+.$$

□

*Remark.* The reason why tame automorphisms are used in this paper is that the author failed to prove the analogue of Lemma 4.3 for general automorphisms.

*Proof of Lemma 4.2.* Let $(A, d') = (T(a_1, \ldots a_n), d')$ be a differential algebra. Assume $\bar{g} \in \mathrm{Aut}_r(A)$. One easily checks that $(A, d')$ is augmented if and only if $(A, (d')^{\bar{g}})$ is. Then, moreover, the linear operators $d'_{(1)}, (d')^{\bar{g}}_{(1)} \colon A_1 \to A_1$ satisfy $(d')^{\bar{g}}_{(1)} = \bar{g}_{(1)} \circ d'_{(1)} \circ \bar{g}_{(1)}^{-1}$, where $\bar{g}_{(1)} = \pi_1 \circ g$, $\pi_1$ being the projection $\bar{A} \to A_1$. This implies $i(A, d') = i(A, (d')^{\bar{g}})$.

Let $g \in \mathrm{Aut}_t(A)$, $g = \bar{g} \circ g_0$, where $g_0 \in \mathrm{Aut}_0(A)$, $\bar{g} \in \mathrm{Aut}_r(A)$. Applying the above argument to $d' = d^{g_0}$, we get the claim of Lemma 4.2.

□

In view of Lemma 4.2, the following assertion is obvious:

**Lemma 4.4** $I((A, d_A) \amalg (B, d_B)) = I(A, d_A) + I(B, d_B)$, *where the sum of two sets is a collection of sums of their elements.*

*Proof of Lemma 4.1.* The invariant $I$ obviously survives tame automorphisms. We need to show that stabilization preserves it. Indeed, by Lemma 4.4, we have $I(S(A, d)) = I(A, d) + I(E, d_E)$. It remains to show that $I(E, d_E) = \{0\}$. Let $g_0 \in \mathrm{Aut}_0(E)$, then $d^{g_0}(e_2) = 0$, and either $d^{g_0}(e_1) = e_2$ or $d^{g_0}(e_1) = e_2 + 1$. Hence $I(E, d_E) = \{0\}$, and the proof is completed. □

*Proof of Theorem 1.1.* We have already checked in Section 3 that the classical invariants of $L_1$ and $L_2$ coincide. We prove now that $I(L_1) = \{3\}$, $I(L_2) = \{1\}$, which would imply that $L_1$ is not isotopic to $L_2$ in the class of Legendrian embeddings.

Consider the differential algebra $(A, d) = (A_{Y_1}, d_{Y_1})$. Let $g \in \mathrm{Aut}_0(A)$ be given by $g(a_i) = a_i + c_i$, $i \in \{1, \ldots, 9\}$. Since $d^g(a_i) = g(d(a_i))$, we have $d^g(a_i) = 0$ when $i \geq 5$ and

$$d_0^g(a_1) = 1 + c_7 + c_5 c_6 c_7,$$
$$d_0^g(a_2) = 1 + c_9 + c_5 c_6 c_9,$$

$$d_0^g(a_3) = 1 + c_7c_8,$$
$$d_0^g(a_4) = 1 + c_8c_9.$$

Assume that $(A, d^g)$ is augmented, that is $d_0^g = 0$. Then $c_7 = c_8 = c_9 = 1$ and $c_5c_6 = 0$. Hence

$$d_{(1)}^g(a_1) = a_7 + c_5a_6 + c_6a_5,$$
$$d_{(1)}^g(a_2) = a_9 + c_5a_6 + c_6a_5,$$
$$d_{(1)}^g(a_3) = a_7 + a_8,$$
$$d_{(1)}^g(a_4) = a_8 + a_9.$$

Since $d_{(1)}^g\left(\sum_{i=1}^4 a_i\right)$ the dimension of the image of $d_{(1)}^g$ equals 3. Hence $I(L_1) = \{3\}$.

Let now $(A, d) = (A_{Y_1}, d_{Y_1})$. Assume that $(A, d^g)$ is augmented. One easily checks that the condition $d_0^g = 0$ implies that $c_7 = c_8 = c_9 = 1$ and $c_5 = 0$. Then

$$d_{(1)}^g(a_1) = a_7 + (c_6 + 1)a_5,$$
$$d_{(1)}^g(a_2) = a_9 + c_6a_5,$$
$$d_{(1)}^g(a_3) = a_7 + a_8,$$
$$d_{(1)}^g(a_4) = a_8 + a_9,$$

and $d_{(1)}^g(a_i) = 0$ when $i \geq 5$. The image of $d_{(1)}^g$ always has dimension 4 and hence $I(L_2) = \{1\}$. □

*Remark.* The construction of the invariant $I$ reflects the following general idea: it is easier to study differential algebras up to (stable tame) automorphisms than their cohomologies up to isomorphisms, and it is easier to study augmented differential free algebras up to isomorphisms preserving the augmentation than general free differential algebras up to all isomorphisms. The invariant $I$ is only the first order obstruction to a stable tame isomorphism. Studying the action of $d^g$ on $\bar{A}/\bar{A}^m$, one can similarly define obstructions (invariants) of order $m$. However, the effect of the stabilization on those higher-order invariants (or at least some of them) is more complicated.

# 5 Proof of Theorem 3.1 and Theorem 3.2

*Proof of Theorem 3.1.* To every double point $a_j$ we associate a positive number $H(a_j)$, called the height of $a_j$. Let $(p_j, q_j, u_j)$, $(p_j, q_j, u'_j)$ be the two points of $L$ projecting onto $a_j$, then $H(a_j) = |u_j - u'_j|$. We start by proving the following important assertion.

**Lemma 5.1** *Let $f \in \widetilde{W}_l(Y)$. Then*

$$\sum_{x \in Q_+} H(f(x)) - \sum_{x \in Q_-} H(f(x)) = \int_{\Pi_k} f^*(dp \wedge dq) > 0, \qquad (2)$$

*where $Q_+$ is the set of positive vertices for $f$, and $Q_-$ the set of negative vertices. In particular, at least one of the vertices of $\Pi_l$ is positive.*

**Corollary 5.2** *Let $(A, d)$, where $A = T(a_1, \ldots, a_n)$, be the differential algebra of a diagram $Y$. Assume that $H(a_1) \leq H(a_2) \leq \cdots \leq H(a_n)$. If $f \in W_k^+(Y, a_{i_1}, \ldots, a_{i_k})$ then $\sum_{j=2}^n H(a_{i_j}) < H(a_{i_1})$. Therefore, $d(a_i) \in T(a_1, \ldots, a_{i-1})$ for every $i \in \{1, \ldots, n\}$.*

This corollary implies that the differential algebra associated with a diagram of Legendrian link is a free model in the sense of [8].

*Proof of Lemma 5.1.* Consider $f \in W_k^+(Y, a_{j_0}, \ldots, a_{j_{k-1}})$. Let $\Gamma_i$ be the piece of $L$ projecting onto $f([x_{i-1}^k, x_i^k]) \subset f(\partial \Pi_k)$, and $\Gamma'_i$ the vertical segment connecting $(p_{j_i}, q_{j_i}, u_{j_i})$ with $(p_{j_i}, q_{j_i}, u'_{j_i})$, where $0 \leq i \leq k-1$. The curves $\Gamma_i, \Gamma'_i$ can be organized into a closed path. Integrating the 1-form $du - p\,dq$ over this path and applying the Stokes' theorem, we get (2). □

Let $f \in W_k^+(Y, a_j, a_{j_0} \ldots, a_{j_{k-1}})$. We have

$$\int_{\Pi_k} f^*(dp \wedge dq) = \sum_{l=1}^m n_l A_l,$$

where $A_l$ is the area of the $l$th component, say $S_l$, of $\mathbf{R}^2 \setminus Y$ and $n_l$ is a nonnegative integer equal to the cardinality of $f^{-1}(z)$, where $z \in S_l$. According to Lemma 5.2, the value of $\sum_{l=1}^m n_l A_l$ does not exceed $H(a_{j_0})$. Hence the coefficients $n_l$ can be chosen in finitely many ways. Thus there are finitely many admissible immersions which send the positive vertex to $a_j$, and the left hand side of (1) contains finitely many terms. □

*Proof of Theorem 3.2.* It suffices to show that $d^2(a_j) = 0$ for every double point $a_j$. Let $N$ be the number of double points. By the definition of $d$, we have
$$d^2(a_j) = \sum_{k \geq 0} \sum_{1 \leq j_1,\ldots,j_k \leq N} C(a_j, a_{j_1}, \ldots, a_{j_k}) \, a_{j_1} \cdots a_{j_k}, \qquad (3)$$
where $C(a_j, a_{j_1}, \ldots, a_{j_k})$ is the cardinality of the set $U(Y, a_j, a_{j_1}, \ldots, a_{j_k})$ of triples $(f', i, f'')$ such that

- $f' \in W_{k'}^+(Y, a_j, a_{j'_1} \ldots, a_{j'_{k'-1}})$;

- $f'' \in W_{k''}^+(Y, a_{j'_i}, a_{j''_1} \ldots, a_{j''_{k''-1}})$;

- $1 \leq i \leq k' - 1$ and 
  $(j_1, \ldots, j_k) = (j'_1, \ldots, j'_{i-1}, j''_1, \ldots, j''_{k''-1}, j'_{i+1}, \ldots, j'_{k'-1})$.

We need to prove that all the numbers $C(a_j, a_{j_1}, \ldots, a_{j_k} \leq N)$ are even.

The proof goes as follows: we introduce some auxiliary space of immersions, $V_k^+(Y, a_{j_0}, \ldots, a_{j_{k-1}})$. Then we construct a map
$$\varphi \colon U(Y, a_{j_0}, \ldots, a_{j_{k-1}}) \to V^+(Y, a_{j_0}, \ldots, a_{j_{k-1}})$$
and two maps
$$\psi_1, \psi_2 \colon V^+(Y, a_{j_0}, \ldots, a_{j_{k-1}}) \to U(Y, a_{j_0}, \ldots, a_{j_{k-1}}),$$
such that the following conditions hold. For all $f \in V^+(Y, a_{j_0}, \ldots, a_{j_{k-1}})$, $\tau \in U(Y, a_{j_0}, \ldots, a_{j_{k-1}})$, we have $f = \varphi(\psi_1(f)) = \varphi(\psi_2(f))$, $\psi_1(f) \neq \psi_2(f)$, and $\tau = \psi_i(\varphi(\tau))$ for $i = 1$ or $i = 2$ (but not both). This would imply that $U(Y, a_{j_0}, \ldots, a_{j_{k-1}})$ has an even number of elements. Indeed, it immediately follows from the above conditions that for every $\tau \in U(Y, a_{j_0}, \ldots, a_{j_{k-1}})$ there exists a unique $f \in V^+(Y, a_{j_0}, \ldots, a_{j_{k-1}})$ and a unique $i \in \{1, 2\}$ such that $\tau = \psi_i(f)$, whence the cardinality of $U(Y, a_{j_0}, \ldots, a_{j_{k-1}})$ is twice that of $V^+(Y, a_{j_0}, \ldots, a_{j_{k-1}})$.

It will be convenient to represent $U(Y, a_j, a_{j_1}, \ldots, a_{j_k})$ as a disjoint union of two subsets, $U^l(Y, a_j, a_{j_1}, \ldots, a_{j_k})$ and $U^r(Y, a_j, a_{j_1}, \ldots, a_{j_k})$ and define the map $\varphi$ separately on each of them.

Let $(f', i, f'') \in U(Y, a_j, a_{j_1}, \ldots, a_{j_k})$. Let $S'$ be a neighbourhood of $x_{j'_i}^{k'}$ in $\Pi_{k'}$ and $S''$ a neighbourhood of $x_0^{k''}$ in $\Pi_{k''}$. Then there are two possible relative positions for the adjacent sectors $S_1 = f'(S')$ and $S_2 = f''(S'')$ near the double point $a_{j'_i} = f'(x^{k'}_{j'_1}) = f''(x^{k''}_0)$ (see Fig. 7).

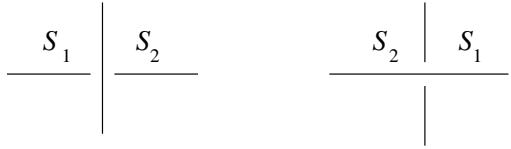

Figure 7:

The triples represented by the left (resp. right) figure belong to $U^l(Y, a_j, a_{j_1}, \ldots, a_{j_k})$ (resp. $U^r(Y, a_j, a_{j_1}, \ldots, a_{j_k})$).

We give the definition of $V^+_k(Y, a_{j_0}, \ldots, a_{j_{k-1}})$. For every $k \geq 1$ fix a (curved) polygon $\Theta_k \subset \mathbf{R}^2$ with $k$ vertices, such that the angle at exactly one of its vertices, say $y^k_0$, exceeds $\pi$.

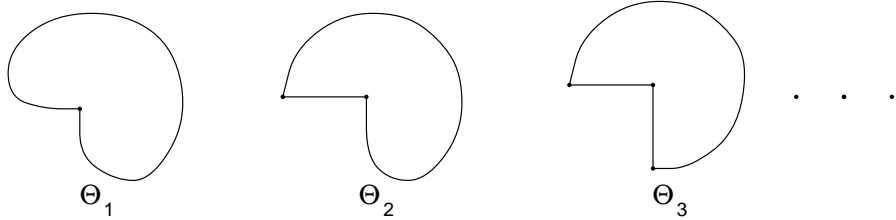

Figure 8:

Let $V_k(Y)$ be the collection of smooth orientation-preserving immersions $f \colon \Theta_k \to \mathbf{R}^2$ such that $f(\partial \Theta_k) \subset Y$ (and hence the images of the vertices are double points of $Y$). Consider the corresponding set of non-parametrized immersions $\widetilde{V}_k(Y) = V_k(Y)/\mathrm{Diff}_+(\Theta_k)$. Let $y^k_0, \ldots, y^k_{k-1}$ be the vertices of $\Theta_k$ numbered counterclockwise. Given $f \in \widetilde{V}_k(Y)$, the definition of positive and negative vertices extends verbatim to the vertices $y^k_i$, $i > 0$. The image of a neighbourhood of $y^k_0$ covers three out of the four sectors in a neighbourhood of $f(y_0)$. We call $y^k_0$ positive if two of them are positive, and negative otherwise. Define the subset $V^+_k(Y)$ of $V_k(Y)$ to consist of such immersions for which exactly one of the vertices $y^k_i$ is positive. The set $V^+_k(Y)$ is a disjoint union of the sets $V^+(Y, a_{j_0}, \ldots, a_{j_{k-1}})$ defined as follows. Let $x^k_s$ be the positive vertex for $f \in V^+_k(Y)$, then $(a_{j_0}, a_{j_1} \ldots, a_{j_{k-1}}) = (f(y^k_s), f(y^k_{s+1}), \ldots, f(y^k_{s-1}))$.

We are going to construct the map $\varphi$. Let $(f', i, f'') \in U^l(Y, a_{j_0}, a_{j_1}, \ldots, a_{j_k})$. Consider parametrized immersions $f'_0 \colon \Pi_{k'} \to \mathbf{R}^2$, $f''_0 \colon \Pi_{k''} \to \mathbf{R}^2$ representing $f', f''$. According to Fig. 7, there exist smooth immersions $r' \colon [0, 1] \to \partial \Pi_{k'}$, $r'' \colon [0, 1] \to \partial \Pi_{k''}$ such that $r'(0) = x_i^{k'}$, $r''(0) = x_0^{k''}$, and $f'_0 \circ r' = f''_0 \circ r''$. Choose the maps $r', r''$ in such a way that their images $r'([0,1]), r''([0,1])$ are the largest possible. Then at least one of the points $r'(1), r''(1)$ is a vertex.

First assume $k'' > 1$. We attach $\Pi_{k'}$ to $\Pi_{k''}$ by identifying $r'(t)$ with $r''(t)$ for every $t \in [0, 1]$. The resulting polygon, which we denote $\Sigma$, is represented in Fig. 9.

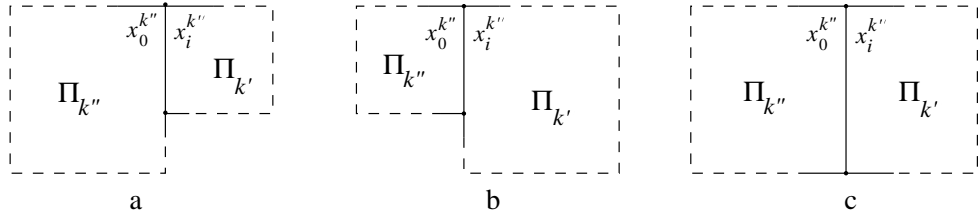

Figure 9:

The version (a) corresponds to the case $r'(1) = x_{i+1}^{k'}$, $r''(1) \neq x_{k''-1}^{k''}$, (b) to the case $r'(1) \neq x_{i+1}^{k'}$, $r''(1) = x_{k''-1}^{k''}$, (c) to the case $r'(1) = x_{i+1}^{k'}$, $r''(1) = x_{k''-1}^{k''}$. The dotted lines are used for the insignificant part of the picture, depending on the actual values of $k'$ and $k''$.

We claim that the case (c) is impossible. Indeed, the result of gluing is a polygon diffeomorphic to $\Theta_{k-2}$. The immersions $f', f''$ can be combined into an immersion $f \in \widetilde{W}_{k-2}(Y)$. All vertices of $\Pi_{k-2}$ are negative, because two of the four vertices which were reduced in the process of gluing had to be positive. This contradicts Lemma 5.1 and hence proves our claim.

In either of the cases (a) and (b), the polygon $\Sigma$ is diffeomorphic to $\Theta_k$, and the immersions $f', f''$ can be combined, in a inique way, into an immersion $f \in \widetilde{V}_k(Y)$. Obviously, exactly one of the vertices of $\Sigma$ (namely that originated from $x_0^{k'}$) is positive, hence $f \in V^+(Y)$. The vertices of $\Sigma$, numbered counterclockwise, are $x_0^{k'}, \ldots, x_{i-1}^{k'}, x_1^{k''}, \ldots, x_{k''-1}^{k''}, x_{i+1}^{k'}, \ldots, x_{k'-1}^{k'}$. Applying $f$ to this sequence, we get the sequence $(a_{j_0}, \ldots, a_{j_{k-1}})$. Thus $f \in V^+(Y, a_{j_0}, \ldots, a_{j_{k-1}})$. Let $\varphi(f', i, f'') = f$.

Assume now that $k'' = 1$. If $r''(1) \neq x_0^1$, we can attach $\Pi_{k'}$ to $\Pi_1$ in the same way as we did in the case $k'' > 1$. The resulting polygon is represented by Fig. 9a, and the definition of $\varphi(f', i, f'')$ is the same as for $k'' > 1$.

Consider the case $r''(1) = x_0^1$, $r'(1) = x_{i+1}^{k'}$. We attach $\Pi_{k'}$ to $\Pi_1$ by identifying $r'(t)$ with $r''(t)$, see Fig. 10a. When $r''(1) = x_0^1$, but $r'(1) \neq x_{i+1}^{k'}$, we have, in addition, to glue together certain pieces of $\Pi_{k'}$. More precisely, let $r_+ : [0,1] \to \partial \Pi_{k'}$, $r_- : [0,1] \to \partial \Pi_{k'}$ be smooth immersions such that $r_+(0) = x_i^{k'}$, $r_-(0) = r'(1)$, and $f_0' \circ r_+ = f_0' \circ r_-$. Choose the maps $r_1, r_2$ to maximize the images $r_+([0,1])$, $r_-([0,1])$. Identify $r_+(t)$ with $r_-(t)$. We have either $r_+(1) = x_{i+1}^{k'}$, $r_-(1) \neq x_{i-1}^{k'}$ (Fig. 10b), or $r_+(1) \neq x_{i+1}^{k'}$, $r_-(1) = x_{i-1}^{k'}$ (Fig. 10c), or $r_+(1) = x_{i+1}^{k'}$, $r_-(1) = x_{i-1}^{k'}$ (Fig. 10d).

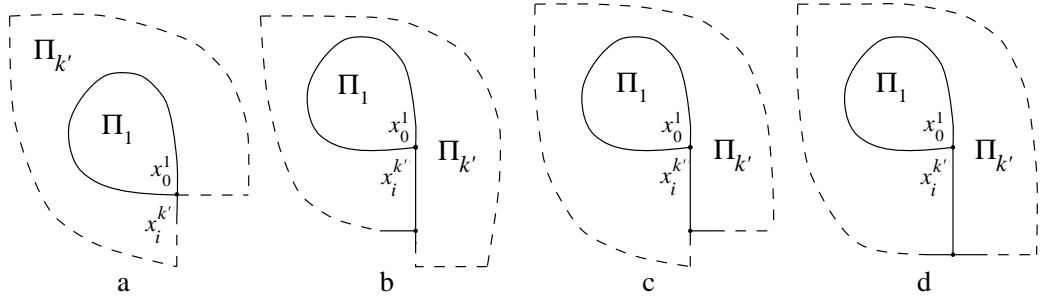

Figure 10:

Actually, the case (d) cannot occur. The reason is the same reason as for Fig. 9c: it would lead to an immersion $f \in \widetilde{W}_{k-2}(Y)$ without positive vertices, which contradicts Lemma 5.1.

In either of the cases (a)–(c), the obtained polygon $\Sigma$ is diffeomorphic to $\Theta_k$. The immersions $f', f''$ can be uniquely combined into an immersion $f \in \widetilde{V}_k(Y)$. Again, exactly one of the vertices of $\Theta_k$ $(x_0^{k'})$ is positive. Hence $f \in V^+(Y)$. The vertices of $\Sigma$, numbered counterclockwise, are $x_0^{k'}, \ldots, x_{i-1}^{k'}, x_0^0, x_{i+1}^{k'}, \ldots, x_{k'-1}^{k'}$, which means that $f \in V^+(Y, a_{j_0}, \ldots, a_{j_{k-1}})$. We define $\varphi(f', i, f'') = f$.

The definition of $\varphi(\tau)$ for $\tau \in U^r(Y, a_{j_0}, \ldots, a_{j_{k-1}})$, is essentialy the same as for $\tau \in U^l(Y, a_{j_0}, \ldots, a_{j_{k-1}})$; actually, they are mirror images of each other.

We are going now to construct the maps $\psi_i$. Let $f \in V^+(Y, a_{j_0}, \ldots, a_{j_{k-1}})$ and fix a parametrized immersion $f_0$ representing $f$. Let immersions $r_1 : [0, \varepsilon] \to \Theta_k$ $r_2 : [0, \varepsilon] \to \Theta_k$ be such that $r_1(0) = r_2(0) = y_0^k$, $r_1([0,1])$ and $r_1([0,1])$ are subsets of $f^{-1}(Y)$, and the segments $r_1([0, \varepsilon])$, $r_2([0, \varepsilon])$ look as represented by the thick line in Fig. 11a,b respectively.

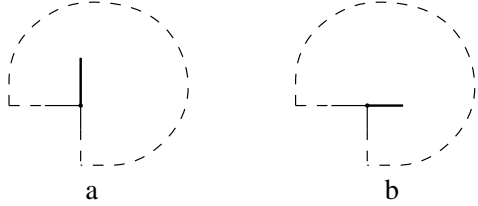

Figure 11:

We extend $r_1, r_2$ to immersions $[0,1] \to \Theta_k$ such that, for each $i \in \{1, 2\}$, $r_i([0,1]) \subset f^{-1}(Y)$, the restriction of $r_i$ to $]0,1[$ is an embedding into the interior of $\Theta_k$, and either $r_i(1) \in \partial \Theta_k$, or $r_i(1) \in r_i([0,1[)$. These immersions are uniquely defined up to reparametrizations. We shall use the immersion $r_1$ ($r_2$) to define $\psi_1(f)$ ($\psi_2(f)$). The possible configurations of the image of $r_1$ are shown in Fig. 12, where the thick line represents $r_1([0,1])$.

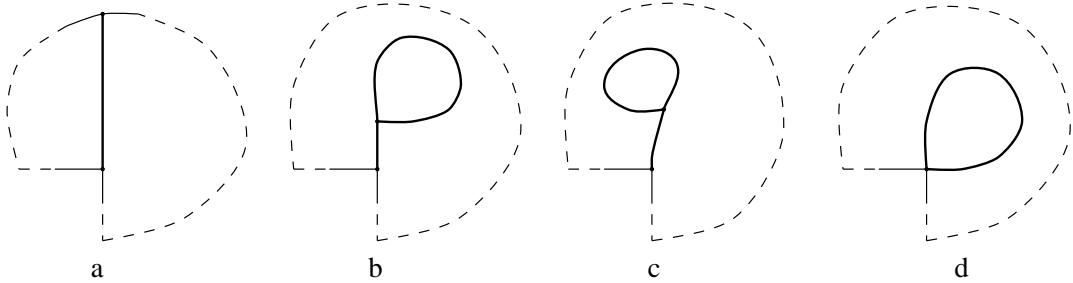

Figure 12:

The set $r_1([0,1])$ divides $\Theta_k$ into two parts, say $\Sigma_1$, $\Sigma_2$, which can be identified with polygons $\Pi_{k'}$, $\Pi_{k''}$, where $k' + k'' = k + 2$. In the case (a), the positive vertex of $\Theta_k$ for the immersion $f$ is a vertex of exactly one of the polygons $\Sigma_1$, $\Sigma_2$. We choose this polygon as $\Pi_{k'}$. In the cases (b)–(d), we choose as $\Pi_{k'}$ the polygon having $k + 1$ vertices.

Restricting $f_0$ to the polygons $\Pi_{k'}$, $\Pi_{k''}$, we get the maps $f_0' \in W_{k'}(Y)$, $f_0'' \in W_{k''}(Y)$. Exactly one vertex of $\Pi_{k'}$ (resp. $\Pi_{k''}$) is positive for $f_0'$ (resp. $f_0''$). This follows from Lemma 5.1 and the fact that the total number of positive vertices equals 2. We can assume that the vertices of the polygons $\Pi_{k'}$, $\Pi_{k''}$ are numbered in such a way that the positive vertices are $x_0^{k'}$, $x_0^{k''}$.

Passing to non-parametrized immersions, we get a pair $f' \in W^+_{k'}(Y, a_{j'_0}, a_{j'_1} \ldots, a_{j'_{k'-1}})$, $f'' \in W^+_{k''}(Y, a_{j''_0}, a_{j''_1} \ldots, a_{j''_{k''-1}})$, Let $i \in \{1, \ldots, k'-1\}$ be such that $x_i^{k'}$ is identified with $x_0^{k''}$ as a point of $\Theta_k$. In the cases (a)–(c), this $i$ is uniquely defined. In the case (d), we choose (say) the smaller of the two possible values of $i$. We define $\psi_1(f) = (f', i, f'')$. It is easy to check that $(f', i, f'') \in U(Y, a_{j_0}, \ldots, a_{j_{k-1}})$.

The definition of $\psi_2$ coincides with that of $\psi_1$ except that we use $r_2$ instead of $r_1$ and, in the case (d) (where, actually, the decompositions of $\Theta_k$ arising from $r_1$ and $r_2$ coincide), we choose the larger of the possible values of $i$.

We need to check that $f = \varphi(\psi_1(f)) = \varphi(\psi_2(f))$ for every $f \in V^+(Y, a_{j_0}, \ldots, a_{j_{k-1}})$. This follows from the easy observation that, decomposing $\Theta_k$ in two parts by means of either $r_1$ or $r_2$ and then gluing it together, we get exactly the same polygon. Similarly, the fact that for every $\tau \in U(Y, a_{j_0}, \ldots, a_{j_{k-1}})$ we have either $\tau = \psi_1(\varphi(\tau))$ or $\tau = \psi_2(\varphi(\tau))$ is equivalent to the following claim, which is easy to check: after glueing together two polygons $\Pi_{k'}, \Pi_{k''}$ into the single polygon $\Theta_k$ as prescribed by the definition of $\varphi$, we can cut $\Theta_k$ using either $r_1$ or $r_2$ to get the same pair of polygons.

We are going to prove that $\psi_1(f) \neq \psi_2(f)$ for all $f$. Suppose, to the contrary, that $\tau = \psi_1(f) = \psi_2(f)$. Since $\tau = \psi_1(f)$, the line of gluing near $y_0^k$ must look as shown on the Fig. 11a. Since $\tau = \psi_2(f)$, it must also look as shown on Fig 11b. Hence the decomposition of $\Theta_k$ has the form represented in Fig. 10d. But then $\psi_1(f) = (f', i, f'') \neq \psi_2(f) = (f', i+1, f'')$.

□

## 6 Proof of Theorem 3.3

We connect $L_0$ with $L_1$ by a generic path $\{L_t\}$ in the space of Legendrian embeddings. Consider the corresponding path $\{Y_t\}$ in the space of diagrams of Legendrian links. In the course of this deformation, the bifurcations (Reidemeister moves) shown in Fig. 13 may take place.

Since the definition of $d$ is purely combinatorial, the differential algebras $(A_{Y_t}, d_{Y_t})$ are naturally isomorphic as long as none of these bifurcations occur. Hence it suffices to show that, $\epsilon$ being small, the differential algebra $(A_{Y_{t'-\varepsilon}}, d_{Y_{t'-\varepsilon}})$ is stable tame isomorphic to $(A_{Y_{t'+\varepsilon}}, d_{Y_{t'+\varepsilon}})$ if at the moment $t = t'$ the diagram goes through one of the bifurcations I–III (note that move IV is the inverse to move III and hence does not require a separate consideration).

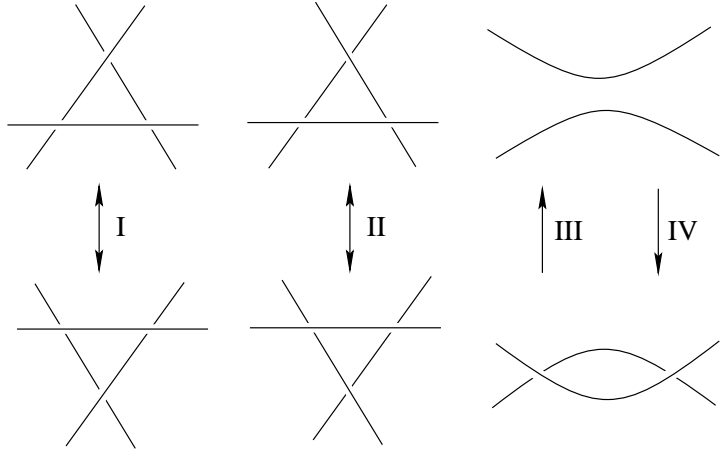

Figure 13:

*Move I.*

Let $(A_{Y_{t'-\varepsilon}}, d_{Y_{t'-\varepsilon}}) = (T(a, b, c, a_1, \ldots, a_l), d_-)$, $(A_{Y_{t'+\varepsilon}}, d_{Y_{t'+\varepsilon}}) = (T(a, b, c, a_1, \ldots, a_l), d_+)$, where $a, b, c$ are as shown in Fig 14, and the points $a_i$ are numbered in the same way for both diagrams. We claim that $d_- = d_+$. Let $f_\pm \in W_3^+(Y_{t'\pm\varepsilon})$ be the immersions whose images are the small triangle with vertices $a, b, c$ These immersions lead to the same term $bc$ in $d_\pm(a)$. It follows from Corollary 5.2 that $H(b) + H(c) < H(a)$.

**Lemma 6.1** *Let $f \in W_k^+(Y_{t'\pm\varepsilon})$, $f \neq f_\pm$. Then neither of the segments $[a, b]$, $[b, c]$, $[a, c]$ is the image of an edge of $\Pi_k$.*

*Proof.* If $[a, b]$ is the image of an edge of $\Pi_k$ then the vertex sent to $b$ is positive and that sent to $a$ is negative. By Corollary 5.2, this is impossible since $H(a) > H(b)$. If $[a, c]$ is the image of an edge of $\Pi_k$ then the vertex sent to $c$ is positive and that sent to $a$ is negative. By Corollary 5.2, this is impossible since $H(a) > H(c)$. If $[b, c]$ is the image of an edge of $\Pi_k$ then the vertices sent to $b$ and $c$ are both positive, which is impossible. □

Assume that $f \in W_k^+(Y_{t'-\varepsilon})$, $f \neq f_\pm$ map one of the vertices of $\Pi_k$ to either $a$, $b$, or $c$. Then $f$ can be smoothly deformed to an immersion $f_+ \in W_k^+(Y_{t'+\varepsilon})$ (two examples are given in Fig 14). Thus we obtain a one-to-one correspondence between $W_k^+(Y_{t'-\varepsilon})$ and $W_k^+(Y_{t'+\varepsilon})$, which leads to the identity $d_- = d_+$.

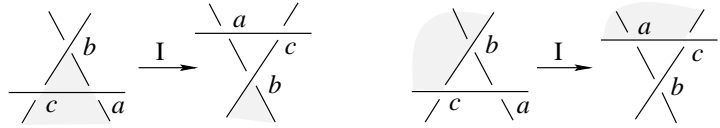

Figure 14:

*Move II.*

Let $(A_{Y_{t'-\varepsilon}}, d_{Y_{t'-\varepsilon}}) = (T(a,b,c,a_1,\ldots,a_l), d_-)$, $(A_{Y_{t'+\varepsilon}}, d_{Y_{t'+\varepsilon}}) = (T(a,b,c,a_1,\ldots,a_l), d_+)$, where $a, b, c$ are as shown in Fig 15, and the points $a_i$ are numbered in the same way for both diagrams. Let $g$ be the elementary automorphism of $T(a,b,c,a_1,\ldots,a_l)$ such that $g$ sends $a$ to $a + bc$ and fixes other generators. We claim that $d_- = d_+^g = g \circ d_+ \circ g^{-1}$ (actually, $g^{-1} = g$) and hence $g$ is a tame isomorphism between $(A_{Y_{t'-\varepsilon}}, d_{Y_{t'-\varepsilon}})$ and $(A_{Y_{t'+\varepsilon}}, d_{Y_{t'+\varepsilon}})$.

Let $f \in W_k^+(Y_{t'\pm\varepsilon})$. The small triangle whose vertices are $a, b, c$ cannot be the image of $f$ since the vertices of $\Pi_3$ sent to $b, c$ would then be positive. It is possible that $f$ sends some of the edges of $\Pi_k$ to the sides of the triangle. If it does not then the immersion extends smoothly through $t = t'$. The immersions $f \in W_k^+(Y_{t'-\varepsilon})$ having one of the fragments shaded in Fig 15a disappear after the bifurcation. Instead, there appear the immersions $f \in W_k^+(Y_{t'+\varepsilon})$ having one of the fragments shaded in Fig 15b. Note that the fragments (1) and (2) (but not (3)) can occur only when the positive vertex is mapped to $a$, and the fragments (3) (but not (1) or (2)) can occur only when the positive vertex is mapped to a double point other then $b$ or $c$.

First we show that $d_-(s) = d_+^g(s) = g(d_+(s))$ when $s \neq a$ is a double point of $Y_{t'-\varepsilon}$. Consider the free algebra $\tilde{A} = T(a^+, a^-, b, c, a_1, \ldots, a_l)$, let $\rho$ be the map $\tilde{A} \to A = T(a, b, c, a_1, \ldots, a_l)$ which sends $a^+$ and $a^-$ to $a$ and fixes other generators. Let $g_+$ (resp. $g_-$) be the automorphism of $\tilde{A}$ which maps $a^+$ to $a^+ + bc$ (resp. $a^-$ to $a^- + bc$) and fixes other generators. We have $g_+ = g_+^{-1}$, $g_- = g_-^{-1}$, $g_+ \circ g_- = g_- \circ g_+ = g_\pm$, and $\rho \circ g_\pm = g$.

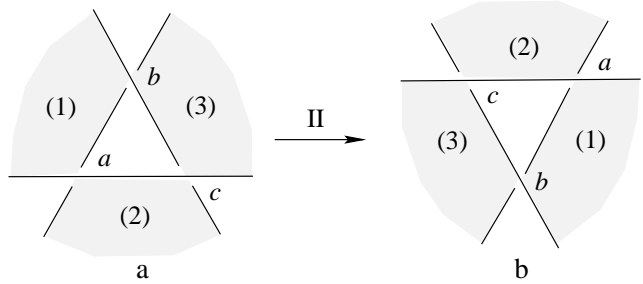

Figure 15:

Let $\xi^+, \xi^- \in \tilde{A}$ be defined as follows. For every $f \in W_k^+(Y_{t'\pm\varepsilon}, s, \ldots)$ we replace each occurrence of $a$ in the corresponding term in $d_\pm(s)$ by $a^-$ if the image of $f$ locally looks as shown in Fig 16a,b and by $a^+$ if it looks as shown in Fig 16c,d. Clearly, we have $\rho(\xi^+) = d_+(s)$, $\rho(\xi^-) = d_-(s)$. Denote by $S_\pm$ the subset of $f \in W_k^+(Y_{t'\pm\varepsilon}, s, \ldots)$ consisting of the immersions which do not contain the fragments marked by (3) in Fig 15. We define $\theta_+, \theta_- \in \tilde{A}$ by summing up the terms in $d_\pm$ corresponding to $f \in S_\pm$. It is easy to check that $\theta_+ = \theta_-$, because these immersions do not contain at all the fragments shown in Fig 15 and hence survive the move.

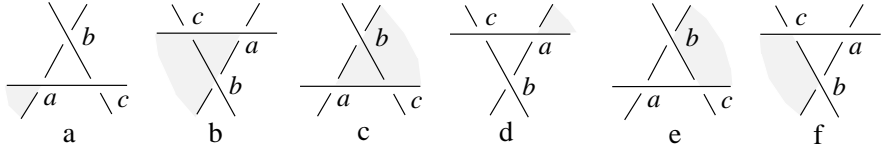

Figure 16:

Each fragment of the kind shown in Fig 16b of an immersion $f \in S_+$ can be replaced by a fragment shown in Fig 16f. Each fragment of the kind shown in Fig 16c of an immersion $f \in S_-$ can be replaced by a fragment shown in Fig 16e. Moreover, all the immersions $f \in W_k^+(Y_{t'\pm\varepsilon}, s, \ldots)$ can be uniquely obtained by applying these procedures. Hence we have $g_+(\theta_+) = \xi^+$, $g_-(\theta_-) = \xi^-$.

We need to prove that $d_-(s) = g(d_+(s))$, or, equivalently,

$$\rho(g_-(\theta_-)) = g(\rho(g_+(\theta_+))).$$

Since $g\circ\rho = \rho\circ g_-\circ g_+$, the right hand side of this equation equals $\rho(g_-(\theta_+))$. Since $\theta_+ = \theta_-$, we get $d_-(s) = g(d_+(s))$.

It remains to show that $d_-(a) = d_+^g(a)$. We can assume that $H(a) > H(b)$ and $H(a) > H(c)$, because $H(a) - H(b) - H(c)$, considered as a function of $t$, vanishes at $t = t'$. Hence, by Corollary 5.2, $a$ is involved in neither $d_+(a)$, nor $d_+(b)$, nor $d_+(c)$. Thus we have $d_+^g(a) = g(d_+(a+bc)) = d_+(a+bc)$. We need to show that $d_-(a) + d_+(a) = d_+(bc)$. Given $i \in \{1,2\}$, let $\lambda_i^+$ (resp. $\lambda_i^-$) be the sum of those terms in $d_+(a)$ (resp. $d_-(a)$) which correspond to the immersions having a fragment marked by $(i)$ in Fig 15a (resp. Fig 15b). We have $d_-(a) + d_+(a) = \lambda_1^+ + \lambda_1^- + \lambda_2^+ + \lambda_2^-$.

For any immersion $f \in W_k^+(Y_{t'+\varepsilon}, c, \ldots)$, we can construct an immersion $f'$ coinciding with $f$ outside a neighbourhood of the small triangle, such that either $f' \in W_k^+(Y_{t'-\varepsilon}, a, b, \ldots)$ or $f' \in W_k^+(Y_{t'+\varepsilon}, a, b, \ldots)$, as shown in Fig 17. This correspondence is one-to-one and hence $\lambda_1^+ + \lambda_1^- = b\, d_+(c)$.

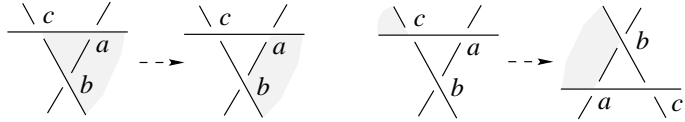

Figure 17:

Similarly, for any immersion $f \in W_k^+(Y_{t'+\varepsilon}, b, \ldots)$, we can construct an immersion $f'$ coinciding with $f$ outside a neighbourhood of the triangle $T$, such that either $f' \in W_k^+(Y_{t'-\varepsilon}, a, \ldots, c)$ or $f' \in W_k^+(Y_{t'+\varepsilon}, a, \ldots, c)$. This correspondence is one-to-one and hence $\lambda_2^+ + \lambda_2^- = d_+(b)c$. Therefore, $d_-(a) + d_+(a) = d_+(b)c + b\, d_+(c)$.

*Move III.*

Let $(A', d') = (A_{Y_{t'-\varepsilon}}, d_{Y_{t'-\varepsilon}})$. Let $a, b$ be the two disappearing double points of $Y_{t'-\varepsilon}$. Perturbing the deformation $Y_t$ and decreasing $\varepsilon$, we can number the double points of $Y_{t'-\varepsilon}$ as $a, b, a_1, \ldots, a_l, b_1, \ldots, b_m$ in such a way that, considering $H(a_i)$, $H(b_j)$ as functions defined on $[t'-\varepsilon, t'+\varepsilon]$ and $H(a)$, $H(b)$ as functions defined on $[t'-\varepsilon, t'[$, we have

$$H(a_l) > \cdots > H(a_1) > H(a) > H(b) > H(b_1) > \cdots > H(b_m).$$

The difference between $H(a)$ and $H(b)$ equals the area of the curved 2-gon $D$ which vanishes at $t = t'$. Let $\bar{f} \in W_2^+(Y_{t'-\varepsilon}, a, b)$ be the immersion whose image is $D$. By Lemma 5.2, the only term in $da$ containing $b$ is that corresponding to $\bar{f}$, and we have $d'(a) = b + v$, where $v \in T(b_1, \ldots, b_m)$. Denote $(A, d) = (A_{Y_{t'+\varepsilon}}, d_{Y_{t'+\varepsilon}})$, where $A = T(a_1, \ldots, a_l, b_1, \ldots, b_m)$. Let $(A_s, d_s) = (A, d) \perp\!\!\!\perp (E, d_E)$.

Consider the isomorphism $\sigma: A' \to A_s$ which sends $a$ to $e_1$, $b$ to $e_2 + v$ and fixes other generators. By definition, the differential algebra $(A', d')$ is tame isomorphic to the differential algebra $(A_s, \hat{d})$, where $\hat{d} = \sigma \circ d' \circ \sigma^{-1}$. We need to prove that there exists a tame isomorphism $g: A_s \to A_s$ such that $d_s = \hat{d}^g$.

Denote by $A_E$ the two-sided ideal in $A_s$ generated by $e_1$ and $e_2$. We have the vector space decomposition $A_s = A \oplus A_E$. Denote $A_{[i]} = T(a_1, \ldots, a_i, b_1, \ldots, b_m, e_1, e_2)$, where $i \in \{0, \ldots, l-1\}$. Let $\gamma: A_s \to A_E$, $\tau: A_s \to A$ be the natural projections.

**Lemma 6.2** a) *The differentials $\hat{d}, d_s$ coincide on $A_{[0]}$.*
b) $\tau \circ (d_s + \hat{d}) = 0$.

Before proving this lemma, we show that it implies the existence of the required isomorphism. Let $g_i \in \mathrm{Aut}_t(A_s)$, $i \in \{1, \ldots, l\}$ be a sequence of elementary automorphisms such that $g_i$ sends $a_i$ to $a_i + q_i$, where $q_i \in A_{[i-1]}$, and fixes all other generators. It gives rise to a sequence $d_{[i]}$ of differentials on $A_s$ such that $d_{[0]} = \hat{d}$, $d_{[i+1]} = d_{[i]}^{g_{i+1}}$. By the definition, all the differential algebras $(A_s, d_{[i]})$ are tame isomorphic. We shall inductively choose the elements $q_i$ in such a way that $d_{[i]}$ would coincide with $d_s$ on $A_{[i]}$. This would yield a tame isomorphism between $(A_s, \hat{d})$ and $(A_s, d_s)$.

We shall need the following lemma concerning stabilizations of differential algebras.

**Lemma 6.3** *There exists a linear map $h: A_s \to A_s$ such that $\gamma = h d_s + d_s h$. Given a subalgebra $A_* \subset A_s$ containing $e_1$ and $e_2$, the map $h$ can be chosen in such a way that $h(A_*) = A_*$.*

This lemma also implies that the projection $\tau$ is chain homotopic to the identity (because $\gamma = \tau + \mathrm{id}$) and hence induces an isomorphism between the homologies of $(A_s, d_s)$ and $(A, d)$.

*Proof of Lemma 6.3.* We use the decomposition $A_s = A \oplus A_E$. For $w \in A$, define $h(w) = 0$. The ideal $A_E$ is spanned, as a linear space, by the elements of the form $x = y e_i z$, where $y \in A$. We define $h(x) = y e_1 z$ if $i = 2$ and $h(x) = 0$ if $i = 1$. The claim is now easy to check. $\square$

Suppose that we have already constructed the automorphisms $g_1, \ldots, g_{j-1}$. Let $g_j$ have the required form, with the value of $q_j \in A_{[j-1]}$ yet to be fixed. We have to choose $q_j$ in such a way that $d_{[j]}(a_j) = d_s(a_j)$.

By the inductive hypothesis, we have $d_{[j-1]}(a_j) = \hat{d}(a_j)$, $d_{[j-1]}(q_j) = d_s(q_j)$. By Corollary 5.2, $d'(a_j) \in T(a, b, a_1, \ldots, a_{j-1}, b_1, \ldots, b_m)$ and hence $\hat{d}(a_j) = \sigma(d'(a_j)) \in A_{[j-1]}$. It also follows from Corollary 5.2 that $d_s(q_j) \in A_{[j-1]}$. Thus $d_{[j-1]}(a_j + q_j) \in A_{[j-1]}$. Since $g_j$ is fixed on $A_{[j-1]}$, this implies

$$d_{[j]}(a_j) = g_j(d_{[j-1]}(a_j + q_j)) = d_{[j-1]}(a_j + q_j) = \hat{d}(a_j) + d_s(q_j). \quad (4)$$

Denote $r = d_s(a_j) + \hat{d}(a_j)$, then $r \in A_{[j-1]}$. By Corollary 5.2, we have $d_s(a_j) \in A_{[j-1]}$. Since $d_s = d_{[j-1]}$ on $A_{[j-1]}$ and $\hat{d}(a_j) = d_{[j-1]}(a_j)$, we get $d_s(r) = d_s^2(a_j) + d_{[j-1]}^2(a_j) = 0$.

Define $q_j = h(r)$. The condition $q_j \in A_{[j-1]}$ is satisfied due to Lemma 6.3. According to (4), it suffices to show that $\hat{d}(a_j) + d_s(q_j) = d_s(a_j)$, which can be written as $r = d_s(q_j)$. By Lemma 6.2, we have $\tau(r) = 0$ and hence $r = \gamma(r)$. It follows from Lemma 6.3 that $\gamma(r) = d_s(q_j) + h(d_s(r))$. Since $d_s(r) = 0$, we get $r = d_s(q_j)$. Hence $(A_s, d_s)$ is tame isomorphic to $(A_s, \hat{d})$.

*Proof of Lemma 6.2.* By the definition of $\sigma$, $d_s(e_1) = \hat{d}(e_1)$ and $d_s(e_2) = \hat{d}(e_2)$. It follows from Corollary 5.2 that the terms $a, b$ are not involved in $d'(b_i)$ for any $i \in \{1, \ldots, m\}$. Hence $\hat{d}(b_i) = d_s(b_i)$. This proves the first claim.

Consider the homomorphism $\zeta: A' \to A$ such that $\zeta(a) = 0$, $\zeta(b) = v$, and $\zeta$ fixes other generators. Then $\tau(\hat{d}(a_j)) = \tau(\sigma(d'(a_j))) = \zeta(d'(a_j))$. Since $d_s(a_j) = d(a_j)$, it remains to show that $\zeta(d'(a_j)) = d(a_j)$ for every $j \in \{1, \ldots, l\}$.

Let $U$ denote the set of immersions $f \in W_k^+(Y_{t'-\varepsilon}, a_j, f(x_1^k), \ldots, f(x_{k-1}^k))$ such that neither of the vertices $x_1^k, \ldots, x_{k-1}^k$ is mapped to $a$. Assume that the vertices $x_{i_1}^k, \ldots, x_{i_p}^k$ are those mapped to $b$ by $f \in U$. Given $f_1, \ldots, f_p$ such that $f_r \in W_{k_r}^+(Y_{t'-\varepsilon}, a, \ldots) \setminus \{\bar{f}\}$, we can construct an immersion $f' \in W_{k'}^+(Y_{t'+\varepsilon}, a_j, \ldots)$, where $k' = k + \sum_{r=1}^p k_r - 2p$, by performing around each of the vertices $x_{i_r}^k$ the operation shown in Fig. 18 (when $p = 0$, we just smoothly deform $f$ into $f'$).

We are going to show that every immersion $f' \in W_{k'}^+(Y_{t'+\varepsilon}, a_j, \ldots)$ can be constructed, in a unique way, by performing this gluing operation. This

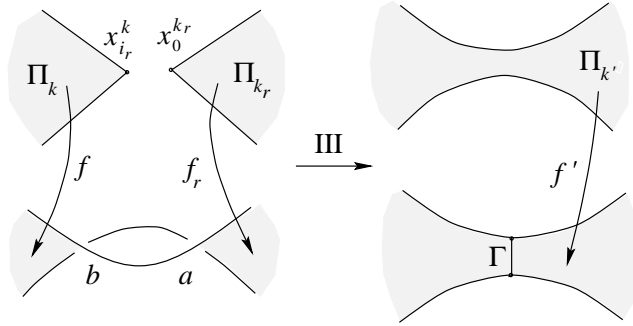

Figure 18:

would imply that

$$\zeta(d'(a_j)) = \sum_{k'\geq 1} \sum_{f'\in \cap W^+_{k'}(Y_{t'+\varepsilon},a_j,f'(x_1^{k'}),...,f'(x_1^{k'}))} f'(x_{k'-1}^{k'}) \cdots f'(x_{k'-1}^{k'}).$$

and hence $\zeta(d'(a_j)) = d_s(a_j)$.

Let $f' \in W^+_{k'}(Y_{t'+\varepsilon}, a_j, \ldots)$. Represent $f'$ by a parametrized immersion $f'_0$. Consider the interval $\Gamma \subset \mathbf{R}^2$ shown in Fig. 18. Let $\Gamma_1, \ldots, \Gamma_n \subset \Pi_{k'}$ be the embedded arcs with ends on $\Pi_{k'}$ such that $f'_0(\Gamma_i) = \Gamma$. Performing around each of the arcs $\Gamma_i$ the operation which is an inverse to that shown in Fig. 18, we construct $n+1$ immersions $f, f_1, \ldots, f_n \in W^+(Y_{t'-\varepsilon})$. Denote by $f \in W^+_k(Y_{t'-\varepsilon})$ the one which inherits from $f'$ the positive vertex $x_0^{k'}$. Let $x_{s_1}^k, \ldots, x_{s_j}^k$ be the vertices sent to $b$ by $f$.

If $j < n$ then there exists, for some $f_i \in W^+_{k_i}(Y_{t'-\varepsilon})$, a negative vertex $x_s^{k_i}$ such that $f_i(x_s^{k_i}) = b$. Since the positive vertex $x_0^{k_i}$ has to be mapped to $a$, we get, by Lemma 5.2, $f_i = \bar{f}$, which is impossible by construction. Therefore, we always have $j = n$.

There exists a canonical numbering of the immersions $f_i \in W^+_{k_i}(Y_{t'-\varepsilon})$ such that, at the moment $t = t'$, $\Pi_{k_i}$ is attached to $\Pi_k$ at $x_{s_i}^k$. Applying the gluing operation to the collection $f, f_1, \ldots, f_n$, we obtain the initial immersion $f'$. It is also easy to check that if $f'$ is a result of gluing together the immersions $f^+, f_1^+, \ldots, f_{n_+}^+$ then actually $(f^+, f_1^+, \ldots, f_{n_+}^+) = (f, f_1, \ldots, f_n)$. This completes the proof of Lemma 6.2.

□

# 7 Some properties of the differential algebra

*7.1. Adding a small kink.* Let $L \subset \mathbf{R}^3$ be an oriented Legendrian knot. Replacing a short interval of its diagram by a kinked curve as shown in Fig. 19, we get two Legendrian knots which we denote $K_-(L)$ and $K_+(L)$.

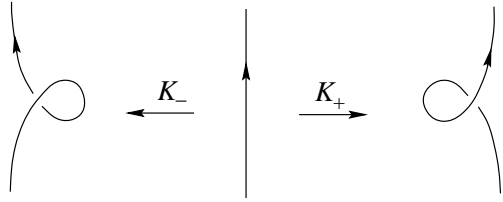

Figure 19:

One easily checks that on the classes of Legendrian isotopy these operations are well-defined (in particular, the resulting Legendrian isotopy class does not depend neither on the interval on which the operation is performed, nor on the size of the kink). Furthermore, this operations commute up to Legendrian isotopy, that is the Legendrian manifolds $K_-(K_+(L))$ and $K_+(K_-(L))$ are Legendrian isotopic.

Denote by $K_\pm^n(L)$ the result of iterating the operation $K_\pm$ $n$ times. It is well-known that if Legendrian knots $L_1$, $L_2$ have the same classical invariants then they become isotopic after adding several kinks: $K_-^{n^-}(K_+^{n^+}(L_1))$ is Legendrian isotopic to $K_-^{n^-}(K_+^{n^+}(L_2))$ (see e. g. [2]). The natural question to ask is how our differential algebras behave under the operations $K_+, K_-$. We call an oriented Legendrian knot reducible if it is Legendrian isotopic to either $K_+(L)$ or $K_-(L)$ for some oriented Legendrian knot $L$, and irreducible otherwise. The following assertion demonstrates that the operations $K_+, K_-$ essentially destroy our invariants.

**Theorem 7.1** *Let $L$ be a reducible Legendrian knot and $(A, d)$ its differential algebra. If $\beta(L)$ is odd then $(A, d)$ is stable tame isomorphic to the differential algebra $(B_{(1)}, d_{B_{(1)}})$, where $B_{(1)} = T(a)$, $d_{B_{(1)}}(a) = 1$. If $\beta(L)$ is even then $(A, d)$ is stable tame isomorphic to the differential algebra $(B_{(2)}, d_{B_{(2)}})$, where $B_{(1)} = T(a, b)$, $d_{B_{(2)}}(a) = 1$ and $d_{B_{(2)}}(b) = 0$.*

Theorem 7.1 implies, in particular, that the cohomology of the differential algebra of a reducible Legendrian knot vanish. Indeed, if $d(z) = 0$ then $d(az) = z$

*Proof of Theorem 7.1.* Assume that $L' = K_+(L'')$ or $L' = K_-(L'')$. Let $(A, d) = (T(a_0, \ldots, a_n), d)$ be the differential algebra associated to the diagram of $L'$, where $a_0$ is the double point of the added the kink. Since $H(a_0)$ equals the area enclosed by the kink, and the latter can be chosen arbitrarily small, it follows from Lemma 5.2 that $d(a_0) = 1$.

Consider the automorphism $g$ which fixes $a_0$ and sends $a_i$ to $a_i + a_0 d(a_i)$, where $i \geq 1$. By Lemma 5.2, $g$ is tame. The differential algebra $(A, d') = (A, d^g)$ is tame isomorphic to $(A, d)$. Since $d(a_0 d(a_i)) = d(a_i)$, we have $d'(a_i) = 0$.

Consider the differential algebra $S^{[n/2]}(B_{(\varepsilon)}, d_{B_{(\varepsilon)}})$, where $\varepsilon$ equals 1 if $n+1$ is odd and 2 if $n+1$ is even. This algebra has $n+1$ generators, since $\beta(L)$ has the same parity as the number of double points, which equals $n+1$. By renaming the generators, it is tame isomorphic to the algebra $(A_s, d_s) = (T(a_0, \ldots, a_n), d_s)$, where $d_s(a_0) = 1$, $d_s(a_{2j-1}) = a_{2j}$ for $j \leq [n/2]$, and $d_s$ vanishes on other generators. Applying the tame automorphism which sends $a_{2j-1}$ to $a_{2j-1} + a_0 a_{2j}$, where $j \leq [n/2]$, and fixes other generators, we obtain an automorphism between the differential algebras $(A_s, d_s)$ and $(A, d')$. Hence there exists a stable tame isomorphism between $(B_{(\varepsilon)}, d_{B_{(\varepsilon)}})$ and $(A, d)$. □

7.2. *Grading.*

The differential algebra $(A, d)$ of a Legendrian link $L$ has a natural grading in which the differential $d$ has degree $-1$. We describe this grading for the simplest case, when $L$ is a knot. In this case, the grading takes values in the group $\mathbf{Z}/m(L)\mathbf{Z}$.

Let $a$ be a double point of the diagram $Y$ of a Legendrian knot $L$. Consider the points $x_+, x_- \in L$ such that $\pi(x_+) = \pi(x_-) = a$ and the point $x_+$ has the larger $u$-coordinate. The points $x_+, x_-$ divide $L$ into two pieces, $C_1$ and $C_2$, which we orient from $x_+$ to $x_-$.

We can assume, without loss of generality, that the intersecting branches are orthogonal at $a$. Then, for $\varepsilon \in \{1, 2\}$, the rotation number of the curve $\pi(C_\varepsilon)$ has the form $N_\varepsilon/2 + 1/4$, where $N_\varepsilon \in \mathbf{Z}$. Clearly, $N_1 - N_2$, up to the sign, is equal to $m(L)$. Hence $N_1$ and $N_2$ represent the same element of the group $\mathbf{Z}/m(L)\mathbf{Z}$, which we define to be the degree of $a$. One can check (perhaps not quite easily) that the differential $d$ has degree $-1$ in this grading.

To extend the definition of a stable tame automorphism to graded al-gebras, we need to modify the definition of stabilization given in Section 2. Instead of the single operation of stabilization, we define a family $S_v$ of them, indexed by $v \in \mathbf{Z}/m(L)\mathbf{Z}$. The differential algebra $S_v(A, d)$ as a non-graded algebra is the same as $S(A, d)$, and we assign degree $v$ to $e_1$ and degree $v-1$ to $e_2$. We call a differential algebra of the form $S_{v_1}(\cdots(S_{v_m}(A, d))\cdots)$ an iterated graded stabilization of the graded differential algebra $(A, d)$. The definition of a graded tame isomorphism can be obtained from that of a tame isomorphism requesting that all the isomorphism involved preserve the grading. We say that two graded differential algebras are stable graded tame isomorphic if some iterated stabilizations of them are graded tame isomorphic.

One easily checks that all algebra homomorphisms which we construct in the proof of Theorem 3.3 respect the grading. Hence the following graded version of Theorem 3.3 holds:

**Theorem 7.2** *Graded differential algebras of Legendrian isotopic Legendrian knots are stable graded tame isomorphic.*

An unpleasant feature of this grading is that usually graded differential algebras of Legendrian knots are not graded isomorphic to augmented graded differential algebras. Thus the machinery which we developed in Section 4 in order to produce the invariants of stable tame isomorphism is incompatible with the grading.

*7.3. Open problems.* We formulate some open problems concerning differential algebras of Legendrian knots.

1. It follows from Theorem 7.1 that reducible Legendrian knots are bad (since their cohomologies vanish). Is it true that all irreducible Legendrian knots are good?

2. Can a reducible Legendrian knot have the same classical invariants as an irreducible one?

3. Consider the involution $F : \mathbf{R}^3 \to \mathbf{R}^3$, $(p, q, u) \mapsto (-p, q, -u)$. It changes the sign of the contact form $\alpha = du - pdq$ and hence maps Legendrian knots to Legendrian knots. Let $L$ be an oriented Legendrian knot. Since $F$ is isotopic to the identity as a linear map, the knots $L$ and $F(L)$ are isotopic as smooth knots. One easily checks that $m(L) = -m(F(L))$ and $\beta(L) = \beta(L)$. The problem suggested by Fuchs and Tabachnikov [2] is: does there exist a Legendrian knot $L$ with $m(L) = 0$ such that $L$ is not Legendrian isotopic to $F(L)$?

Let $(A_L, d_L)$ be the differential algebra of a Legendrian knot $L$. Then the differential algebra $(A_{F(L)}, d_{F(L)})$ associated to $F(L)$ can be obtained by applying to $(A, d)$ the automorphism $F_*$ of $A = T(a_1, \ldots, a_n)$ which reverses the order of the factors: $F_*(a_{i_1} \cdots a_{i_k}) = a_{i_k} \cdots a_{i_1}$. Our question is: is it possible for $(A_{F(L)}, d_{F(L)})$ not to be stable tame isomorphic to $(A_L, d_L)$?

e-mail: yuri@chekanov.mccme.rssi.ru